\documentclass[11pt]{article}

\usepackage{latexsym}
\usepackage{amssymb}

\newtheorem{thm}{Theorem}[section]
\newtheorem{lem}{Lemma}[section]
\newtheorem{prop}{Proposition}[section]
\newtheorem{cor}{Corollary}[section]

\begin{document}
{

\begin{center}
\Large\bf
On a J-polar decomposition of a bounded operator and matrix representations of J-symmetric,
J-skew-symmetric operators.
\end{center}
\begin{center}
\bf S.M. Zagorodnyuk
\end{center}

\noindent
{\large\bf Introduction.}

\noindent
Complex symmetric, skew-symmetric and orthogonal matrices are classical objects of
the finite-dimensional linear analysis~\cite{Cit_1000_Gantmacher}. In particular,
the canonical spectral forms are known for them.
Certainly, they have a more complicated structures as for Hermitian matrices.
However, in a certain sense complex symmetric matrices are more universal.
Namely, an {\it arbitrary} square complex matrix is similar to a symmetric matrix.
If one introduces a J-form and write conditions for a symmetric, skew-symmetric and orthogonal
matrix (continued by zeros to the right and to the bottom to obtain a semi-infinite matrix) in its terms, one
arrives to the well-known J-symmetric, J-skew-symmetric and J-isometric
operators.

A general definition of a J-symmetric operator was given by I.M.~Glazman in his paper~\cite{Cit_2000_GlazmanStDAN}.
A study of these operators had been continued in papers of N.A.~Zhyhar and A.~Galindo
(see the references in a monograph~\cite{Cit_3000_GlazmanKniga}).
Later, an investigation of these operators had been performed by A.D.~Makarova, L.A.~Kamerina,
V.P.~Li, T.B.~Kalinina, A.N.~Kochubey, B.G.~Mironov (a seria of papers by these authors
appeared in 70-th, 80-th of the 20-th century in Ulyanovskiy sbornik "Funkcionalniy analiz"),
L.M.~Rayh, E.R.~Tsekanovskiy and others.
Most of these papers were devoted to the questions of extensions of J-symmetric operators to
J-self-adjoint operators and to a description of all such extensions.
At the present time, J-self-adjoint operators are studied by S.R.~Garcia, M.~Putinar,
E.~Prodan (see the paper~\cite{Cit_4000_GarciaPutinarCS2} and References therein).

A definition of a bounded J-skew-symmetric operator was given by Sh.~Asadi and I.E.~Lutsenko in the
paper~\cite{Cit_5000_AsadiLutsenko}.
A general definition appeared in a paper of T.B.~Kalinina~\cite{Cit_6000_Kallinina1981},
she continued to study these operators in papers~\cite{Cit_7000_Kallinina1982},~\cite{Cit_8000_Kallinina1983}.
J-symmetric and J-skew-symmetric operators also appeared in a book~\cite{Cit_9000_GohbergKreinTVO1967}
in a study of Volterra operators context.

In papers of L.A.~Kamerina J-isometric and quasi-unitary operators and a notion of
quasi-unitary equivalence
were introduced~\cite{Cit_10000_Kamerina1975},\cite{Cit_11000_Kamerina1987}.

Consider a separable Hilbert space $H$. Recall that a conjugation (involution) in
$H$ is an operator $J$, defined on the whole $H$ è and satisfying the following
properties~\cite{Cit_12000_AkhGlazmanT1TLO1977},\cite{Cit_13000_StoneLO}
\begin{equation}
\label{f1_0}
J^2 = E,\qquad (Jx,Jy) = \overline{(x,y)},\qquad x,y\in H,
\end{equation}
where $E$ is the identity operator in $H$, and $(\cdot,\cdot)$ is a scalar product in $H$.
For each conjugation there exists an orthonormal basis $\mathcal F = \{ f_k \}_{k\in{\mathbb Z}_+}$ in $H$
such that
\begin{equation}
\label{f1_1}
J x = \sum_{k=0}^\infty \overline{x_k} f_k,
\qquad x=\sum_{k=0}^\infty x_k f_k\in H.
\end{equation}
This basis is not uniquely determined, it is determined up to a unitary transformation which
commutes with $J$ (J-real).
An arbitrary such a basis $\mathcal F$ we shall call {\bf corresponding} to the involution $J$.
Define the following linear with respect to the both arguments functional (J-form):
\begin{equation}
\label{f1_2}
[x,y]_J := (x,Jy),\qquad x,y\in H.
\end{equation}
A linear operator $A$ in $H$ is said to be J-symmetric, if
\begin{equation}
\label{f1_3}
[Ax,y]_J = [x,Ay]_J,\qquad x,y\in D(A),
\end{equation}
and is said to be J-skew-symmetric if
\begin{equation}
\label{f1_4}
[Ax,y]_J = -[x,Ay]_J,\qquad x,y\in D(A).
\end{equation}
If the following condition is true:
\begin{equation}
\label{f1_5}
[Ax,Ay]_J = [x,y]_J,\qquad x,y\in D(A),
\end{equation}
then the operator is said to be J-isometric.

Let the domain of $A$ is dense in $H$. The operator $A$ is said to be J-self-adjoint if
\begin{equation}
\label{f1_6}
A=JA^*J,
\end{equation}
and is said to be J-skew-self-adjoint if
\begin{equation}
\label{f1_7}
A=-JA^*J.
\end{equation}
If
\begin{equation}
\label{f1_8}
A^{-1} = JA^*J,
\end{equation}
then the operator $A$ we shall call a $J$-unitary.
Notice that the operator $A^T = JA^*J$ in \cite{Cit_12000_AkhGlazmanT1TLO1977}
was called {\bf transposed} (later, in some papers it was also called J-adjoint,
but we shall use the latter word for the operator $\widetilde A = JAJ$).

For non-densely defined operators, one can also introduce a notion of J-symmetric and
J-skew-symmetric linear relations, see, e.g.,~\cite{Cit_14000_Mironov1987}.

Let $A$ be a linear bounded operator in $H$.
In this case, conditions~(\ref{f1_3}),(\ref{f1_4}),
(\ref{f1_5}) mean that the matrix of the operator in an arbitrary basis $\mathcal F$,
which is coorresponding to $J$,
will be symmetric, skew-symmetric and orthogonal, respectively.
This remark and some properties of the J-form allow to obtain some simple properties of
eigenvalues and eigenvectors of such matrices.

In this work we obtain a J-polar decomposition for bounded operators (under some conditions).
This decomposition is analogous to the polar decomposition of a bounded operator and to
the J-polar decomposition in J-spaces~\cite{Cit_15000_KrShmJpolarn}.
Also, we  obtain other decompositions which are analogous to decompositions for
finite-dimensional matrices in~\cite{Cit_1000_Gantmacher}.
A possibility of the matrix representation for J-symmetric and J-skew-symmetric operators and its
properties are studied.
A structure of the following null set
$H_{J;0} = \{ x\in H: [x,x]_J = 0 \}$,
is studied, as well.

\noindent
{\bf Notations. }
As usual, we denote by ${\mathbb R}, {\mathbb C}, {\mathbb N}, {\mathbb Z}, {\mathbb Z}_+, {\mathbb R}^2$
the sets of real numbers, complex numbers, positive integers, non-negative integers and
the real plane, respectively.
Everywhere in this paper, all Hilbert spaces are assumed to be separable,
$(\cdot,\cdot)$ and $\| \cdot \|$ denote the scalar product and the norm in a Hilbert space, respectively.

For a set $M$ in a Hilbert space $H$, by $\overline M$
we mean a closure of $M$ in the norm $\| \cdot \|$.
For $\{ x_k \}_{k\in{\mathbb Z}_+},\ x_k\in H,$ we write
$\mathop{\rm Lin}\nolimits \{ x_k \}_{k\in{\mathbb Z}_+} := \{ y\in H:\ y=\sum_{j=0}^n \alpha_j x_j,\
\alpha_j\in{\mathbb C},\ n\in{\mathbb Z}_+ \}$;
$\mathop{\rm span}\nolimits \{ x_k \}_{k\in{\mathbb Z}_+} := \overline{\mathop{\rm Lin}\nolimits \{ x_k \}_{k\in{\mathbb Z}_+}}$.

The identity operator in a Hilbert space $H$ is denoted by $E$.
For an arbitrary linear operator $A$ in $H$, the operators
$A^*, \overline A, A^{-1}$ mean its adjoint operator, its closure and its  inverse (if they exist).
By $D(A)$ and $R(A)$ we mean the domain and the range of the operator $A$, and by
$\mathop{\rm Ker}\nolimits A$ we mean the kernel of the operator $A$. By $\sigma(A),\ \rho(A)$
we denote the spectrum of $A$ and the resolvent set of $A$, respectively.
The  resolvent function of $A$ we denote by $R_\lambda (A)$, $\lambda\in\rho (A)$.
Also, we denote $\Delta_A(\lambda) = (A-\lambda E) D(A)$.
The norm of a bounded operator $A$ is denoted by $\| A \|$.

By $l_2$ we denote the space of complex sequences
$x=(x_0,x_1,x_2,...)^T$, $x_k\in{\mathbb C}$, $k\in{\mathbb Z}_+$, with a finite norm
$\| x \| = \left( \sum_{k=0}^\infty |x_k|^2 \right)^{\frac{1}{2}}$ (the superscript $T$ stands for the
transposition).
\section{Properties of eigenvalues and eigenvectors.}
We shall begin with some simple properties of J-symmetric, J-skew-symmetric and J-orthogonal
operators which, in particular, lead to some new properties of finite-dimensional
complex symmetric, skew-symmetric and orthogonal matrices.
Let $J$ be a conjugation in a Hilbert space $H$.

Vectors $x$ and $y$ are said to be {\bf J-orthogonal}, if $[x,y]_J=0$.
The following proposition is true (concerning statement (i) of the Proposition see.~Theorem~2 in
a paper~\cite[p.86]{Cit_16000_Li1974}).
\begin{prop}
\label{t2_1}
Let $A$ be a J-symmetric operator in a Hilbert space $H$. The following statements are true:

\noindent
$\mathrm{(i)}$ Eigenvectors of the operator $A$ which correspond to different eigenvalues are
J-orthogonal;

\noindent
$\mathrm{(ii)}$ If vectors $x$ and $Jx$, $x\in D(A)$, are eigenvectors of the operator $A$,
then they correspond to the same eigenvalue.
\end{prop}
{\bf Proof. }
In fact, we can write
$$\lambda_x [x,y]_J = [Ax,y]_J = [x,Ay]_J = \lambda_y [x,y]_J, $$
and therefore
\begin{equation}
\label{f2_1}
(\lambda_x - \lambda_y) [x,y]_J = 0.
\end{equation}
Suppose that $x,\overline{x}:=Jx \in D(A)$ are eigenvectors of the operator $A$,
which correspond to eigenvalues
$\lambda_x$  and $\lambda_{\overline x}$, respectively.
Write~(\ref{f2_1}) with $y=\overline{x},\lambda_y=\lambda_{\overline x}$:
$$ (\lambda_x - \lambda_{\overline{x}}) [x,\overline{x}]_J=0. $$
Since $[x,\overline{x}]_J = \| x \|^2 >0$, we get $\lambda_x = \lambda_{\overline x}$.
$\Box$

Define the following set:
\begin{equation}
\label{f2_0}
H_{J;0} := \{ x\in H:\ [x,x]_J = 0 \}.
\end{equation}
In a similar to the latter proof manner the validity of the following two propositions is established.
\begin{prop}
\label{t2_2}
If $A$ is a J-skew-symmetric operator in a Hilbert space $H$, then the following is true:

\noindent
$\mathrm{(i)}$ Eigenvectors of the operator $A$, which correspond to non-zero eigenvalues,
belong to the set $H_{J;0}$;

\noindent
$\mathrm{(ii)}$ If $\lambda_x,\lambda_y$ are eigenvalues of the operator $A$ such that
$\lambda_x \not= -\lambda_y$, then the corresponding to them eigenvectors are J-orthogonal;

\noindent
$\mathrm{(iii)}$ Suppose that $x,\overline{x}:=Jx \in D(A)$ are eigenvectors of the operator
$A$, corresponding to the eigenvalues
$\lambda_x$  and $\lambda_{\overline x}$, respectively. Then $\lambda_x = -\lambda_{\overline x}$.
\end{prop}
\begin{prop}
\label{t2_3}
Let $A$ be a J-isometric operator in a Hilbert space $H$. Then the following statements are true:

\noindent
$\mathrm{(i)}$ Eigenvectors of the operator $A$, which correspond to different from $\pm 1$ eigenvalues
belong to the set $H_{J;0}$;

\noindent
$\mathrm{(ii)}$ If $\lambda_x,\lambda_y$ are eigenvalues of the operator $A$ such that
$\lambda_x \not= \frac{1}{\lambda_y}$, then the corresponding to them eigenvectors are J-orthogonal;

\noindent
$\mathrm{(iii)}$ Suppose that $x,\overline{x}:=Jx \in D(A)$ are eigenvectors
of the operator $A$, corresponding to the eigenvalues
$\lambda_x$ and $\lambda_{\overline x}$, respectively. Then $\lambda_x =
\frac{1}{\lambda_{\overline x}}$.
\end{prop}
It is interesting to notice that in the finite-dimensional case the point  $0$ for a skew-symmetric
matrix and points $\pm 1$ for an orthogonal matrix are distinguished in a special manner in
the spectrum, as well.

In the case of a unitary space $U^n$ with a dimension $n$, $n\in{\mathbb Z}_+$, in an analogous manner,
a conjugation $J$, a J-form, and J-orthogonality are defined. So,
the latter statements are true for complex symmetric, skew-symmetric and orthogonal matrices.

\noindent
{\bf Example 1.1. }
Consider a numerical matrix
$A = \left( \begin{array}{cc} 1 & i\\
i & 0\end{array}\right)$.
Its eigenvalues are
$\lambda_1 = \frac{1}{2} + \frac{\sqrt{3}}{2} i$,
$\lambda_2 = \frac{1}{2} - \frac{\sqrt{3}}{2} i$,
and the corresponding to them normed eigenvectors are
$f_1 = \frac{1}{2\sqrt{2}}\left( \begin{array}{cc}  \sqrt{3} - i\\
2\end{array}\right)$,
$f_2 = \frac{1}{2\sqrt{2}}\left( \begin{array}{cc}  -\sqrt{3} - i\\
2\end{array}\right)$.
Vectors $f_1,f_2$ are not orthogonal. However, they are J-orthogonal.

Let $J$ be a conjugation in a Hilbert space $H$ and $A$ be a bounded linear operator in $H$.
The norm of $A$, as it can be easily seen from the properties of the involution, can be calculated
by the following formula
\begin{equation}
\label{f2_2}
\| A \| = \sup_{x,y\in H:\ \|x\|=\|y\|=1} |[A x,y]_{J}|.
\end{equation}
The following statement is true:
\begin{prop}
\label{p2_0}
If $A$ is a bounded J-symmetric operator in a Hilbert space $H$, then its norm can be calculated as
\begin{equation}
\label{f2_3}
\| A \| = \sup_{x\in H:\ \|x\| = 1} | [Ax,x]_{J} |.
\end{equation}
\end{prop}
{\bf Proof. }
Consider an operator $A$ such as in the statement of the Proposition.
Set
$C := \sup_{x\in H:\ \|x\| = 1} | [Ax,x]_J |$.
For arbitrary elements $x,y\in H:\ x\not=\pm y$ we can write
$$[A(x+y),x+y]_J - [A(x-y),x-y]_J = 4 [Ax,y]_J; $$
$$| [Ax,y]_J | \leq \frac{1}{4} \left( |[A(x+y),x+y]_J| + |[A(x-y),x-y]_J|
\right) = $$
$$= \frac{1}{4} ( \left| \left[ A(\frac{x+y}{\| x+y \|}),\frac{x+y}{\| x+y \|}\right]_J
\right| \| x+y \|^2 +
\left| \left[ A(\frac{x-y}{\| x-y \|}),\frac{x-y}{\| x-y \|} \right]_J\right| * $$
\begin{equation}
\label{f2_3_1}
* \| x-y \|^2 ) \leq \frac{1}{4} C (\|x+y\|^2 + \|x-y\|^2) = \frac{1}{2} C (\| x\|^2 + \| y\|^2).
\end{equation}
Thus, by using~(\ref{f2_2}) and~(\ref{f2_3_1}) we get
$$
\| A\| = \sup_{x,y\in H:\ \|x\|=\|y\|=1} | [Ax,y]_J | \leq C.
$$
On the other hand, we can write
$$
C = \sup_{x,y\in H:\ \|x\|=1} | [Ax,x]_J | \leq
\sup_{x,y\in H:\ \|x\|=\|y\|=1} | [Ax,y]_J | =\| A\|.
$$
Therefore $C=\|A\|$. $\Box$

For a J-skew-symmetric operator $A$, its norm can not be calculated by the formula~(\ref{f2_3}).
Moreover, the following characteristic property of J-skew-symmetric operators is true.
\begin{prop}
\label{p2_1}
A linear operator $A$
in a Hilbert space $H$ is J-skew-symmetric if and only if the following equality is true
\begin{equation}
\label{f2_4}
[Ax,x]_J = 0,\qquad x\in D(A).
\end{equation}
\end{prop}
{\bf Proof. }
We first notice that from the properties of an involution it follows that
$[x,y]_{J} = [y,x]_{J},\quad x,y\in H$.
Let us check the necessity. From relation~(\ref{f1_4}) it follows that
$$[Ax,x]_{J} = - [x,Ax]_{J} =
-[Ax,x]_{J}, $$
and therefore~(\ref{f2_4}) holds true.

\noindent
Let us check the sufficiency.
By using~(\ref{f2_4}) we write
$$0 = [A(x+y),x+y]_{J} = [Ax,x]_{J} + [Ax,y]_{J}
+ [Ay,x]_{J} + [Ay,y]_{J} =$$
$$=[Ax,y]_{J} + [Ay,x]_{J},\qquad x,y\in D(A). $$
From this relation we obtain that $[Ax,y]_{J} = - [Ay,x]_{J}
= -[x,Ay]_{J}$.
$\Box$

Let $J$ be a conjugation in a Hilbert space $H$ and $A$ be an arbitrary linear operator
in $H$.
The operator $\widetilde A := \widetilde{(A)}_J := JAJ$  we shall call
{\bf J-adjoint} to the operator $A$. We first note that $\widetilde{\widetilde A}
= A$ and the following easy to check lemma is true.
\begin{lem}
\label{l3_1}
For a linear operator $A$ in a Hilbert space $H$,
equalities $\overline{D(A)}=H$ and $\overline{D(\widetilde A)}=H$
are true or false simultaneously. The same can be said about equalities
$\overline{R(A)}=H$ and $\overline{R(\widetilde A)}=H$.
\end{lem}

An operation of the construction of the J-adjoint operator commutes with
basic operations on operators. Let us formulate the necessary for us
properties as propositions.
\begin{prop}
\label{p3_1}
Let $A$ be a linear operator in a Hilbert space $H$ such that
$\overline{D(A)}=H$ and $J$ be a conjugation in $H$. Then the following relation is true
\begin{equation}
\label{f3_1}
\widetilde{A^*} = (\widetilde A)^*.
\end{equation}
\end{prop}
{\bf Proof. }
Choose an arbitrary element $g\in D((\widetilde A)^*)$. On one hand, it is true
$$(\widetilde Ax, g) = (x,(\widetilde A)^* g) = (JJx,JJ(\widetilde A)^*g) =
\overline{(Jx,J(\widetilde A)^*g)} = $$
$$= (J(\widetilde A)^* g, Jx),\qquad x\in D(\widetilde A). $$
On the other hand, wee can write
$$(\widetilde Ax, g) = (JAJx,JJg) = \overline{(AJx,Jg)} = (Jg, AJx),\qquad
x\in D(\widetilde A). $$
Comparing right hand sides we obtain that
$$(AJx, Jg) = (Jx, J(\widetilde A)^*g), $$
and therefore $Jg\in D(A^*)$, $A^* Jg = J(\widetilde A)^*g$. Multiplying by $J$ both sides
of the latterr equality we get $\widetilde{A^*} g = (\widetilde A)^*g$.
Therefore
\begin{equation}
\label{f3_3}
(\widetilde A)^* \subseteq \widetilde{A^*}.
\end{equation}
In order to obtain the inverse inclusion one should write the inclusion~(\ref{f3_3})
with the operator $\widetilde A$, and then to take J-adjoint operators for the both sides
(the inclusion under the last operation will stay true).
$\Box$

\begin{prop}
\label{p3_2}
Let $A$ be a linear operator in a Hilbert space $H$ and $J$ be a conjugation in $H$.
Suppose that operators $A$ and $\widetilde A$ admit closures. Then
the following equality is true
\begin{equation}
\label{f3_4}
\widetilde{\overline{A}} = \overline{\widetilde A}.
\end{equation}
\end{prop}
{\bf Proof. }
Choose an arbitrary element $g\in D(\overline{\widetilde A})$.
Then there exists a sequence $x_n \in D(\widetilde A),n\in{\mathbb Z}_+,$ such that
$x_n\rightarrow x$, $\widetilde A x_n = JAJ x_n \rightarrow \overline{\widetilde A} x$
when $n\rightarrow\infty$.
By continuity of the operator $J$ from this it follows that
$$Jx_n \rightarrow Jx,\quad A J x_n\rightarrow J \overline{\widetilde A} x. $$
Consequently, we obtain, that
$Jx\in D(\overline{A})$ è $\overline{A} Jx = J \overline{\widetilde A} x$.
Therefore $x\in D(\widetilde{\overline{A}})$ and
$\widetilde{\overline{A}} x = \overline{\widetilde A} x$. From this relation we conclude that
\begin{equation}
\label{f3_5}
\overline{\widetilde A} \subseteq \widetilde{\overline{A}}.
\end{equation}
In order to obtain the inverse inclusion, we write the inclusion~(\ref{f3_5})
for the operator $\widetilde A$, and then to take  J-adjoint operators for the both sides.
$\Box$

\begin{prop}
\label{p3_3}
Let $A$ be a linear invertible operator in a Hilbert space $H$ and $J$ be a conjugation in $H$.
Then the operator $\widetilde A$ is also invertible and the following equality is true
\begin{equation}
\label{f3_6}
\widetilde{A^{-1}} = (\widetilde A)^{-1}.
\end{equation}
\end{prop}
{\bf Proof. }
Since
$\widetilde{A^{-1}} \widetilde{A} = E|_{D(\widetilde{A})}$,
and $D(\widetilde{A^{-1}}) = J D(A^{-1}) = JR(A) = R(\widetilde{A})$,
the operator $\widetilde A$ is invertible and relation~(\ref{f3_6}) is true.
$\Box$

Notice that a condition of a J-symmetric operator~(\ref{f1_3}) with the help of J-adjoint operator
will be written as follows:
\begin{equation}
\label{f3_7}
(Ax,y) = (x, \widetilde A y),\qquad x\in D(A),\ y\in D(\widetilde A).
\end{equation}
Conditions of a J-skew-symmetric operator~(\ref{f1_4}) and J-isometric operator(\ref{f1_5}) will be
written as
\begin{equation}
\label{f3_8}
(Ax,y) = -(x, \widetilde A y),\qquad x\in D(A),\ y\in D(\widetilde A),
\end{equation}
and
\begin{equation}
\label{f3_9}
(Ax,\widetilde A y) = (x,y),\qquad x\in D(A),\ y\in D(\widetilde A),
\end{equation}
respectively.

Now we shall assume that the operator $A$ is densely defined in $H$.
Notice that in this case from condition~(\ref{f3_9}) it follows that the operator
$A$ is invertible. In fact, equality $Ax=0$ implies the equality
$(x,y)=0$ on a dense in $H$ set $D(\widetilde A)$. Thus,
a densely defined J-isometric operator is always invertible.

Note that in the case of a densely defined operator $A$, conditions~(\ref{f3_7}),(\ref{f3_8}),(\ref{f3_9})
are equivalent to the following conditions
\begin{equation}
\label{f3_13}
A \subseteq (\widetilde A)^*,
\end{equation}
\begin{equation}
\label{f3_14}
A \subseteq -(\widetilde A)^*,
\end{equation}
and
\begin{equation}
\label{f3_15}
A^{-1} \subseteq (\widetilde A)^*,
\end{equation}
respectively.
From these relations, in particular, it immediately follows that densely defined J-symmetric and
J-skew-symmetric operators admit closures. As it is seen from relations~(\ref{f1_3}),(\ref{f1_4}),
their closures will also be J-symmetric or J-skew-symmetric operators, respectively.
For a densely defined J-isometric operator one can only state that
its inverse operator admits a closure.
However, from relation~(\ref{f1_5}) it is easily seen that the inverse operator
to a J-isometric is also J-isometric.
Consequently, if the range of the original J-isometric operator
(the domain of the inverse operator) is also dense,
then it admits a closure.
In this case, also from the relation~(\ref{f1_5}),
it is seen that this closure will be a J-isometric operator.

Note that the operation of the construction of a J-adjoint operator
does not change the defined above by us types of operators.
Namely, the following proposition is true:
\begin{prop}
\label{p3_4}
Let $A$ be a linear operator in a Hilbert space $H$ and $J$ be a conjugation in $H$.
If the operator $A$ is J-symmetric, J-skew-symmetric or J-isometric,
then the same is the operator $\widetilde A = J A J$, as well.
\end{prop}
{\bf Proof. }
The statement about a J-symmetric (J-skew-symmetric, J-isometric)
operator follows from relation~(\ref{f3_7}) ((\ref{f3_8}), (\ref{f3_9})),
respectively, taking into account that $A=\widetilde{\widetilde A}$.
$\Box$

For an element $x\in H$ and a set $M\subseteq H$ we write
$x \perp_J M$, if $x\perp_J y$, forr all $y\in M$.
For a set $M\subseteq H$ we denote $M^\perp_J
= \{ x\in H:\ x\perp_J y,\ y\in M \}$.

It is known that the residual spectrum of a J-self-adjoint operator is empty. It follows
from the theorem below.
\begin{thm} (\cite[Theorem 4, p.87]{Cit_16000_Li1974})
\label{t3_1}
Let $A$ -be a J-self-adjoint operator in a Hilbert space $H$.
A complex number $\lambda$ is an eigenvalue of $A$ if and only if
\begin{equation}
\label{f3_16}
\overline{\Delta_A(\lambda)} \not= H.
\end{equation}
In this case, $(\Delta_A(\lambda))^\perp_J$ will be an eigen-subspace which corresponds to $\lambda$.
\end{thm}
We shall obtain analogous results for J-skew-symmetric and J-isometric operators.
The following theorem is true:
\begin{thm}
\label{t3_2}
Let $A$ be a J-skew-self-adjoint operator in a Hilbert space $H$.
A complex value $\lambda$ is an eigenvalue of $A$ if and only if
\begin{equation}
\label{f3_21}
\overline{\Delta_A(-\lambda)} \not= H.
\end{equation}
In this case, $(\Delta_A(-\lambda))^\perp_J$ will be an eigen-subspace which corresponds to $\lambda$.
\end{thm}
{\bf Proof. }
{\it Necessity. }
Let $x$ be an arbitrary eigenvector of the operator $A$ which corresponds to an eigenvalue $\lambda$. Since $A$,
in particular,  is skew-symmetric, then we can write
for an arbitrary $y\in D(A)$
\begin{equation}
\label{f3_17}
0 = [(A-\lambda E)x, y]_J = -[x,(A+\lambda E) y]_J.
\end{equation}
Therefore $x\perp_J \Delta_A(-\lambda)$ and by the continuity of $[\cdot,\cdot]_J$
we get
\begin{equation}
\label{f3_18}
x\perp_J \overline{\Delta_A(-\lambda)}.
\end{equation}
Since $[x,Jx]=\| x\|^2>0$, then $Jx \notin \overline{\Delta_A(-\lambda)}$ and
therefore $\overline{\Delta_A(\lambda)} \not= H$.

\noindent
{\it Sufficiency. }
Suppose that equality (\ref{f3_21}) is true.
Then there exists $0\not=y\in H$ such that
\begin{equation}
\label{f3_19}
(z,y)=0,\qquad z\in\overline{\Delta_A(-\lambda)}.
\end{equation}
Therefore
$((A+\lambda  E)x,y)=0$, and from this relation we get
$(Ax,y)=(x,\overline{(-\lambda)} y)$, $x\in D(A)$.
Thus, we have $y\in D(A^*)$ and
\begin{equation}
\label{f3_20}
A^* y = -\overline{\lambda} y.
\end{equation}
But since $A$ is J-skew-self-adjoint, then $A^*=-\widetilde A$, and we obtain
$$\widetilde A y = \overline{\lambda} y. $$
From this relation it follows that
$Jy \not= 0$ is an eigenvector of the operator $A$ with an eigenvalue $\lambda$.

Let us show that the following set
\begin{equation}
\label{f3_16_1}
V(\lambda) := (\Delta_A(-\lambda))^\perp_J \backslash \{ 0 \},
\end{equation}
is a set of eigenvectors of the operator $A$,
corresponding to a eigenvalue $\lambda$.
Denote the latter set by $S(\lambda)$.
By the proven property~(\ref{f3_18}), the inclusion $S(\lambda) \subseteq V(\lambda)$ is true.
On the other hand, if $x\in V(\lambda)$, then for $y:= Jx$ relation~(\ref{f3_19}) is true.
Repeating arguments which follow after this formula
we come to a conclusion that $x$ is an eigenvector of the operator $A$,
corresponding to
$\lambda$. Thus, the inverse inclusion is also true.

Finally, since $A=(\widetilde A)^*$, then $A$ is closed. Therefore
$(\Delta_A(-\lambda))^\perp_J$ is an eigen-subspace of the operator
$A$, which corresponds to $\lambda$.
$\Box$
\begin{cor}
\label{c3_1}
The point $0$ can not belong to the residual spectrum of a J-skew-self-adjoint operator.
\end{cor}
In an analogous manner, the following result for J-unitary operators is established.
\begin{thm}
\label{t3_3}
Let $A$ be a J-unitary operator in a Hilbert space $H$.
A complex number $\lambda$ is an eigenvalue of $A$ if and only if
\begin{equation}
\label{f3_22}
\overline{\Delta_A\left( \frac{1}{\lambda} \right)} \not= H.
\end{equation}
In this case, $(\Delta_A(\frac{1}{\lambda}))^\perp_J$ is an eigen-subspace,
which corresponds to $\lambda$.
\end{thm}
\begin{cor}
\label{c3_2}
Points $\pm 1$ can not belong to the residual spectrum of a J-unitary operator.
\end{cor}

From relations~(\ref{f3_7}),(\ref{f3_8}) it is seen that a defined in the whole $H$ J-symmetric (J-skew-symmetric)
operator is a bounded J-self-adjoint (respectively J-skew-swlf-adjoint) operator.
The following statements are also true.
\begin{prop} (\cite[Theorem 1, p.85-86]{Cit_16000_Li1974},\cite[Theorem 3, p.69]{Cit_6000_Kallinina1981})
\label{p3_5}
Let $A$ be  a linear densely defined operator in a Hilbert space $H$, which is
J-symmetric (J-skew-symmetric).
Suppose that $R(A)=H$. Then the operator $A$
is a J-self-adjoint (respectively J-skew-self-adjoint) operator.
\end{prop}
\begin{prop}
\label{p3_6}
Let $A$ be a linear densely defined operator in a Hilbert space $H$, which is
J-symmetric (J-skew-symmetric).
Suppose that $\overline{R(A)}=H$. Then the operator $A$ is invertible and the operator $A^{-1}$
is also  a J-symmetric (respectively J-skew-symmetric) operator.
\end{prop}
{\bf Proof. }
In a view of analogous considerations, we shall check the validity of this Proposition
only for the case of a J-skew-symmetric operator $A$.
Notice that
$\mathop{\rm Ker}\nolimits A^* = H\ominus \overline{R(A)} = \{ 0 \}$.
Thus, the operator $A^*$ is invertible. Since $A$ is J-skew-symmetric, the following inclusion is true
$\widetilde A\subseteq -A^*$ and therefore $\widetilde A$ is invertible, as well.
By Proposition~\ref{p3_3} we conclude that the operator $A$ has an inverse operator.
From  the inclusion $\widetilde A\subseteq -A^*$ it follows the following inclusion
\begin{equation}
\label{f3_25}
(\widetilde A)^{-1} \subseteq -(A^*)^{-1}.
\end{equation}
Notice that $\overline{D(A^{-1})} = \overline{R(A)} = H$. Thus, we can state that
 $(A^*)^{-1} = (A^{-1})^*$.
Using this equality and using Proposition~\ref{p3_3},
from relation~(\ref{f3_25}) we obtain the following inclusion
$$\widetilde{A^{-1}} \subseteq -(A^{-1})^*. $$
And this means that the operator $A^{-1}$ is J-skew-symmetric.
$\Box$
\begin{prop}
\label{p3_7}
Let $A$ be a J-self-adjoint (J-skew-self-adjoint) operator in a Hilbert space $H$.
Suppose that $\overline{R(A)}=H$. Then the operator $A$ is invertible and the operator
$A^{-1}$
is also J-self-adjoint (respectively J-skew-self-adjoint) operator.
\end{prop}
{\bf Proof. }
In a view of analogous considerations, we shall give the proof only for the case
of J-self-adjoint operator $A$.
By Proposition~\ref{p3_6} the operator $A$ is invertible.
By Proposition~\ref{p3_3} the operator $\widetilde A$ is invertible, as well.
From Lemma~\ref{l3_1} it follows that $\overline{R(\widetilde A)} = H$ è $\overline{D(\widetilde A)}=H$.
Thus, we have $\overline{D((\widetilde A)^{-1})}=H$.
Consequeently, the following equality is true
$((\widetilde A)^*)^{-1} = ((\widetilde A)^{-1})^*$.
Since the operator $A$ is J-self-adjoint, the last equality can be written as
$A^{-1} = ((\widetilde A)^{-1})^*$.
Using Proposition~\ref{p3_3}, we obtain the following equality
$A^{-1} = (\widetilde{A^{-1}})^*$, which shows that the operator $A^{-1}$
is J-self-adjoint.
$\Box$

\section{A J-polar decomposition of bounded operators.}
We shall extend in the case of J-symmetric, J-skew-symmetric and J-isometric operators
a seria of properties of finite-dimensional complex symmetric, skew-symmetric and orthogonal
matrices (see~\cite{Cit_1000_Gantmacher}).

The following lemma is true:
\begin{lem}
\label{l5_1}
Let $A$ be a bounded self-adjoint and J-isometric operator in a Hilbert space $H$. Then
the operator $A$ admits the following representation:
\begin{equation}
\label{f5_1}
A = I e^{iK},
\end{equation}
where $I$ is a bounded self-adjoint J-real involutory ($I^2=E$) operator in $H$,
and $K$ is a commuting with $I$ bounded skew-self-adjoint J-real operator in $H$.

If additionally it is known that the operator $A$ is positive, $A\geq 0$, then
one can choose $I=E$.
\end{lem}
{\bf Proof. }
Consider an operator $A$ such as in the statement of the Lemma.
Since the operator $A$ is J-isometric and bounded, then from~(\ref{f1_5})
we obtain
$A^* J A = J$, or $A^* \widetilde A = E$. Since $A$ is self-adjoint, then
\begin{equation}
\label{f5_2}
A \widetilde A = E.
\end{equation}
For the operator $A$ we can write the following representation
\begin{equation}
\label{f5_3}
A = S + iT,
\end{equation}
where $S = \frac{1}{2}(A+\widetilde A),\ T = \frac{1}{2i}(A-\widetilde A)$.
By this, operators $S$ and $T$ are J-real, the operator $S$ is self-adjoint and J-self-adjoint,
and the operator $T$ is skew-self-adjoint and J-skew-self-adjoint.
Since $\widetilde A = S - iT$, then from relation~(\ref{f5_2}) we get
$$E = A \widetilde A = (S+iT)(S-iT) = S^2 + T^2 + i(TS-ST). $$
From this relation it follows that operators $T$ and $S$ commute and
\begin{equation}
\label{f5_4}
S^2 + T^2 = E.
\end{equation}
Since operators $S$ and $iT$ are commuting bounded self-adjoint operators, then they admit
spectral representations
\begin{equation}
\label{f5_5}
S = \int_L \lambda dE_\lambda,\qquad iT = \int_L z dF_z,
\end{equation}
where $E_\lambda,\ F_z$ are commuting resolutions of unity of the operators, and
$L=(l_1,l_2]$, $l_1,l_2\in{\mathbb R},$ is a finite interval of the real line which contains
the spectra of operators.
From equality~(\ref{f5_4}), by using spectral resolutions we get
\begin{equation}
\label{f5_7}
\int_L \int_L (\lambda^2 - z^2 - 1) dE_\lambda dF_z = 0,
\end{equation}
where the integral means
a limit in the norm of $H$ of the corresponding Riemann-Stieltjes type sums (in the  plane).

A point $(\lambda_0,z_0) \in {\mathbb R}^2$ we call {\bf a point of increase} for the measure $dE_\lambda dF_z$, if
for an arbitrary number $\varepsilon>0$,
there exists an element $x\in H$ such that
\begin{equation}
\label{f5_8}
( E_{ \lambda_0 + \varepsilon } - E_{ \lambda_0 - \varepsilon } ) ( F_{ z_0 + \varepsilon } - F_{ z_0 - \varepsilon } ) x \not= 0,
\end{equation}
or, equivalently,
\begin{equation}
\label{f5_9}
( ( E_{ \lambda_0 + \varepsilon } - E_{ \lambda_0 - \varepsilon } ) ( F_{ z_0 + \varepsilon } - F_{ z_0 - \varepsilon } ) x, x) > 0.
\end{equation}
For an arbitrary point of increase
$(\lambda_0,z_0) \in {\mathbb R}^2$ of the measure $dE_\lambda dF_z$ it  is true
\begin{equation}
\label{f5_9_1}
\lambda_0^2 - z_0^2 - 1 = 0.
\end{equation}
In fact,
if the latter equality is not true for a point of increase
$u_0 = (\lambda_0,z_0) \in {\mathbb R}^2$, then
$| \lambda^2 - z^2 - 1 | \geq a,\ a>0,$ in a neighborhood
$U = U(\lambda_0,z_0;\varepsilon) = \{ (\lambda,z)\in{\mathbb R}^2:\ \lambda_0-\varepsilon < \lambda \leq \lambda_0 + \varepsilon,\
z_0 - \varepsilon < z \leq z_0 + \varepsilon \}$, $\varepsilon>0$,
of the point $u_0$.
For this number $\varepsilon$, there exists an element $x\in H$ such that~(\ref{f5_9}) is true.
But
$$0= \| \int_L \int_L (\lambda^2 - z^2 -1) dE_\lambda dF_z x \|^2
= \int_L \int_L |\lambda^2 - z^2 -1 |^2 (dE_\lambda dF_z x,x) \geq $$
$$\geq \int\int_{U} |\lambda^2 - z^2 -1 |^2 (dE_\lambda dF_z x,x) \geq
a^2 ((E_{\lambda_0 + \varepsilon} - E_{\lambda_0 - \varepsilon})(E_{z_0 + \varepsilon} - E_{z_0 - \varepsilon}) x, x) > 0. $$
If two continuous functions
$\varphi (\lambda,z)$ è $\widehat\varphi (\lambda,z)$
on $L^2 = \{ (\lambda,z)\in{\mathbb R}^2:\ \lambda,z\in L \}$ coincide in the points of increase of the measure
$dE_\lambda dF_z$, then
\begin{equation}
\label{f5_10}
\int_{L} \int_{L} \varphi (\lambda,z) dE_\lambda dF_z = \int_{L} \int_{L} \widehat \varphi (\lambda,z) dE_\lambda dF_z.
\end{equation}
In fact,
$$ \| \int_{L} \int_{L} (\varphi (\lambda,z) - \widehat\varphi (\lambda,z)) dE_\lambda dF_z x \|^2
= \int_{L} \int_{L} | \varphi (\lambda,z) - \widehat\varphi (\lambda,z) |^2 (dE_\lambda dF_z x,x), $$
and it remains to notice that $(dE_\lambda dF_z x,x)$
is a positive measure on $L^2$, and
the function under the integral
is equal to zero in all points of increase of this measure.

Consider a set $\Gamma \subset {\mathbb R}^2$, which consists of points $(\lambda,z) \in {\mathbb R}^2$, such that
\begin{equation}
\label{f5_11}
\lambda^2 - z^2 - 1 = 0.
\end{equation}
From~(\ref{f5_11}) it follows that for all poiints of the set $\Gamma$ it is true
$|\lambda| = \sqrt{1+z^2}$ (where we mean the arithmetic value of the root). Hence, for all points of $\Gamma$
\begin{equation}
\label{f5_12}
\lambda = \mathrm{sgn}(\lambda) \sqrt{1+z^2},
\end{equation}
where
\begin{equation}
\label{f5_13}
\mathrm{sgn}(\lambda) = \left\{ \begin{array}{cc} 1, & \lambda > 0,\\
-1, & \lambda \leq 0 \end{array}\right. .
\end{equation}
By the identity $z=\mathop{\rm sh}\nolimits\mathop{\rm arcsh}\nolimits z$, the equality~(\ref{f5_12})
can be rewritten in the following form
\begin{equation}
\label{f5_14}
\lambda = \mathrm{sgn}(\lambda) \sqrt{
\mathop{\rm ch}\nolimits^2( \mathop{\rm arcsh}\nolimits z )} = \mathrm{sgn}(\lambda)
\mathop{\rm ch}\nolimits ( \mathop{\rm arcsh}\nolimits z ),
\end{equation}
in a view of positivity of the hyperbolic cosine function.
By this representation we can write
$$ A = S+iT = \int_{L} \int_{L} (\lambda +z) dE_\lambda dF_z =
\int_{L} \int_{L} ( \mathrm{sgn}(\lambda) \mathop{\rm ch}\nolimits ( \mathop{\rm arcsh}\nolimits z ) + z) dE_\lambda dF_z = $$
\begin{equation}
\label{f5_15}
= \int_{L} \int_{ L_+ } e^{\mathop{\rm arcsh}\nolimits z} dE_\lambda dF_z +
\int_{L} \int_{ L_- } (-e^{-\mathop{\rm arcsh}\nolimits z}) dE_\lambda dF_z,
\end{equation}
where $L_+ = (0,\infty) \cap L,\ L_- = (-\infty,0] \cap L$.

Define the following operator
\begin{equation}
\label{f5_16}
V = \int_{L} \int_{L} \mathrm{sgn}(\lambda) \mathop{\rm arcsh}\nolimits z dE_\lambda dF_z =
\int_{L} \mathrm{sgn}(\lambda) dE_\lambda  \int_{L} \mathop{\rm arcsh}\nolimits z dF_z.
\end{equation}
The operator $V$ is bounded self-adjoint and J-imaginary.
In fact, since the operator $S$ is
J-real, then its resolution of unity $E_\lambda$ commutes with $J$ (see~\cite{Cit_13000_StoneLO}). Therefore
the operator
\begin{equation}
\label{f5_17}
I := \int_{L} \mathrm{sgn}(\lambda) dE_\lambda,
\end{equation}
is a bounded J-real self-adjoint involutory operator.
On the other hand, $\mathop{\rm arcsh}\nolimits (iT) = \sum_{k=0}^\infty a_{2k+1} (iT)^{2k+1},$ $a_{2k+1}\in{\mathbb R}$,
is a J-imaginary, as a limit of J-imaginary operators (here the convergence is
understood in the norm of $H$).

From relations~(\ref{f5_15}),(\ref{f5_16}),(\ref{f5_17}) we conclude that
$$A= Ie^V. $$
Set $K = -iV$, and we obtain the required representation~(\ref{f5_1}).

If it is additionally known that the operator $A$ is positive, $A\geq 0$, then
$$I = A e^{-V} = ( e^{-\frac{V}{2}} )^* A e^{-\frac{V}{2}},$$
is positive, as well.
Therefore $I$ is a positive square root of $E$. By the uniqueness of such a root
we conclude that $I=E$.
$\Box$

The following theorem is true.
\begin{thm}
\label{t5_1}
Let $A$ be a bounded J-unitary operator in a Hilbert space $H$.
The operator $A$ admits the following representation:
\begin{equation}
\label{f5_18}
A = R e^{iK},
\end{equation}
where $R$ is J-real unitary operator in $H$,
and $K$ is a bounded J-real skew-self-adjoint operator in $H$.
\end{thm}
{\bf Proof. }
Consider an operator $A$ such as in the statement of the Theorem.
Suppose that representation~(\ref{f5_18}) is true. Then
$$ A^* A = e^{iK} R^* R e^{iK} = e^{2iK}. $$
Now we shall do not assume an existence of representation~(\ref{f5_18}) and notice
that the operator $G := A^* A$ is
positive self-adjoint and J-unitary.
In fact, since the operator $A$ is bounded by assumption and J-unitary, then
$A^*$ is also bounded and J-unitary. A product of bounded J-unitary operators
is a bounded J-unitary operator, this is verified directly by the definition.
By Lemma~\ref{l5_1} we find a bounded
J-real skew-self-adjoint operator $K$ such that
\begin{equation}
\label{f5_19}
G =  e^{2iK}.
\end{equation}
Now set by definition
\begin{equation}
\label{f5_20}
R = A e^{-iK}.
\end{equation}
By equality~(\ref{f5_19}) we can write
$$ R^* R = e^{-iK} A^*A e^{-iK} = E, $$
and, hence, the operator $R$ is unitary.
Now notice that
$$ J e^{-iK} J = J (\cos(iK) - i \sin(iK) ) J = \cos(iK) + i \sin(iK), $$
since the operator $iK$ is J-real and therefore its resolution of unity commutes with $J$.
Consequently, we have
\begin{equation}
\label{f5_21}
J e^{-iK} J = e^{iK} = (e^{-iK})^{-1},
\end{equation}
and the operator $e^{-iK}$ is J-unitary.
By~(\ref{f5_20}),(\ref{f5_21}) and using that the operator $A$
is J-unitary we conclude that
$$ R^{-1} = e^{iK} A^{-1} = J e^{-iK} J J A^* J = J (A e^{-iK})^* J = \widetilde{(R^*)}, $$
and therefore the operator $R$ is J-unitary.
Then $R^{-1} = R^* = JR^*J$, and therefore $R^*$ is a J-real operator.
Using  matrix representations of operators $R^*$ and $R$ in an arbitrary basis,
which corresponds to the involution $J$,
we conclude that the operator $R$ is J-real.
$\Box$

\begin{lem}
\label{l5_2}
Let $A$ be a J-self-adjoint and unitary operator in a Hilbert space $H$.
The operator $A$ admits the following representation:
\begin{equation}
\label{f5_22}
A = e^{iS},
\end{equation}
where $S$ is a bounded J-real self-adjoint operator in $H$.
\end{lem}
{\bf Proof. }
Consider an operator $A$ such as in the statement of the Lemma.
For the J-self-adjoint operator $A$ it is true $A^*=\widetilde A$, and we can write
the following representation
\begin{equation}
\label{f5_23}
A = S + iT,
\end{equation}
where $S = \frac{1}{2}(A+\widetilde A) = \frac{1}{2}(A+ A^*),\
T = \frac{1}{2i}(A-\widetilde A) = \frac{1}{2i}(A- A^*)$.
Here operators $S$ and $T$ are J-real and self-adjoint.
Since the operator $A$ is unitary, then
$$E = A^* A = (S-iT)(S+iT) = S^2 + T^2 + i(ST-TS). $$
From this relation it follows that operators $T$ and $S$ commute and
\begin{equation}
\label{f5_24}
S^2 + T^2 = E.
\end{equation}
Since operators $S$ and $T$ are commuting bounded self-adjoint operators, then they admit
the following spectral resolutions
\begin{equation}
\label{f5_25}
S = \int_L \lambda dE_\lambda,\qquad T = \int_L z dF_z,
\end{equation}
where $E_\lambda,\ F_z$ are commuting resolutions of unity of operators, and
$L=(l_1,l_2]$, $l_1,l_2\in{\mathbb R},$ is a finite interval of the real line, which contains
the spectra of operators. Moreover, since operators $S$ and $T$
are J-real, then their resolutions of unity commute with $J$.
By equality~(\ref{f5_24}) and using spectral resolutions we get
\begin{equation}
\label{f5_26}
\int_L \int_L (\lambda^2 + z^2 - 1) dE_\lambda dF_z = 0,
\end{equation}
where the integral means a limit in the norm of $H$ of the corresponding Riemann-Stieltjes type sums.
Thus, in all points of increase of the measure $dE_\lambda dF_z$ the following relation is true
\begin{equation}
\label{f5_27}
\lambda^2 + z^2 - 1 = 0.
\end{equation}
A circle~(\ref{f5_27}) in the plane ${\mathbb R}^2$ we denote by $\Gamma$.
For all points of the circle $\Gamma$ it is true
$|z| = \sqrt{1-\lambda^2}$ (where we mean the arithmetic value of the root).
Therefore for all points of $\Gamma$
\begin{equation}
\label{f5_28}
z = \mathrm{sgn}(z) \sqrt{1-\lambda^2},
\end{equation}
where $\mathrm{sgn}(\cdot)$ is from~(\ref{f5_13}).
By the identity $\lambda=\cos\arccos \lambda,$ $\lambda\in [-1,1],$ the equality~(\ref{f5_28})
can be rewritten in the following form
\begin{equation}
\label{f5_29}
z = \mathrm{sgn}(z) \sqrt{\sin^2( \arccos \lambda )} = \mathrm{sgn}(z) \sin ( \arccos \lambda ),
\end{equation}
where we have used the positivity of sine function on $[0,\pi]$.
By this representation we can write
$$ A = S+iT = \int_{L} \int_{L} (\lambda +iz) dE_\lambda dF_z =
\int_{L} \int_{L} ( \cos\arccos \lambda + i \mathrm{sgn}(z) * $$
\begin{equation}
\label{f5_30}
* \sin (\arccos \lambda) ) dE_\lambda dF_z = \int_{ L_+ } \int_{ L } e^{i\arccos \lambda} dE_\lambda dF_z +
\int_{ L_- } \int_{ L } e^{-i\arccos \lambda} dE_\lambda dF_z,
\end{equation}
where $L_+ = (0,\infty) \cap L,\ L_- = (-\infty,0] \cap L$.
Define the following operator
\begin{equation}
\label{f5_31}
S := \int_{ L } \int_{ L } \mathrm{sgn}(z) \arccos \lambda dE_\lambda dF_z =
\int_{ L } \mathrm{sgn}(z) dF_z \int_{ L } \arccos \lambda dE_\lambda.
\end{equation}
It is obvious that $S$ is a J-real self-adjoint operator.
From relation~(\ref{f5_30}) it is seen that~(\ref{f5_22}) is true.
$\Box$

Using the proven lemma we shall establish the following theorem.
\begin{thm}
\label{t5_2}
Let $A$ be a unitary operator in a Hilbert space $H$.
The operator $A$ admits the following representation:
\begin{equation}
\label{f5_32}
A = R e^{iS},
\end{equation}
where $R$ is J-real unitary operator in $H$,
and $S$ is a bounded J-real self-adjoint operator in $H$.
\end{thm}
{\bf Proof. }
Consider an operator $A$ such as in the statement of the Theorem.
Suppose that representation~(\ref{f5_32}) is true. Then
$A^* = e^{-iS} R^*$ è
$$ \widetilde{A^*} = \widetilde{ e^{-iS} } \widetilde{ R^* } =
J( \cos S - i \sin S )J \widetilde{ R^* } = ( \cos S + i \sin S ) \widetilde{ R^* } = e^{iS} R^*, $$
since $S$ and $R$ are J-real.
Since $R$ is unitary, we can write
\begin{equation}
\label{f5_33}
\widetilde{A^*} A = e^{iS} R^* R e^{iS} = e^{2iS}.
\end{equation}
Now we shall not suppose that representation~(\ref{f5_32}) holds true.
Since the operator $A$ is unitary, then
operators $A^{-1} = A^*$, $JA^*J$ and $G := \widetilde{A^*} A$ are unitary, as well.
The operator $G$ is J-self-adjoint since
$G^* = A^* \widetilde{A} = J \widetilde{A^*} A J = \widetilde{G}$.
Applying to this operator Lemma~\ref{l5_2} we find J-real self-adjoint operator $S$
such that
\begin{equation}
\label{f5_34}
G = e^{2iS}.
\end{equation}
Now we set by definition
\begin{equation}
\label{f5_35}
R = A e^{-iS}.
\end{equation}
The operator $R$ is unitary as a product of two unitary operators.
Then we can write
$\widetilde{R^*} = J e^{iS} A^* J = e^{-iS} \widetilde{A^*}$,
and therefore
$$ \widetilde{R^*} R = e^{-iS} \widetilde{A^*} A e^{-iS} = e^{-iS} G e^{-iS} = E. $$
Since the range of a unitary operator $R$ is the whole $H$, then by the latter equality
we get $\widetilde{R^*} = R^{-1}$. Thus, the operator $R$ is J-unitary.
Since the operator $R$ is unitary and J-unitary, it is J-real.
From~(\ref{f5_35}) it follows the representation~(\ref{f5_32}).
$\Box$

Let $A$ be a linear bounded operator in a Hilbert space $H$ and $J$ be a conjugation in $H$.
It is easy to verify that operators $A^T A = JA^* J A$, $A A^T = A JA^* J$
are bounded J-self-adjoint operators.
If $A^T A = A A^T$, then the operator $A$ we shall call {\bf J-normal}.
It is clear that bounded J-self-adjoint, J-skew-self-adjoint and J-unitary operators are J-normal.
The following theorem is true:
\begin{thm}
\label{t5_3}
Let $A$ be a linear bounded operator in a Hilbert space $H$ and $0\notin \sigma(A)$.
Let $J$ be a conjugation in $H$.
Suppose that the spectrum of the operator $A A^T$ has an empty intersection with
a radial ray
$L_\varphi =
\{ z\in{\mathbb C}:\ z=xe^{i\varphi},\ x\geq 0 \}$ ($\varphi\in [0,2\pi)$) in the complex plane.
Then the operator $A$ admits a representation
\begin{equation}
\label{f5_36}
A = S U,
\end{equation}
where $S$ is a bounded J-self-adjoint operator in $H$, and $U$ is a bounded J-unitary operator in $H$.
Here
\begin{equation}
\label{f5_37}
S = \sqrt{AA^T},
\end{equation}
where the square root is understood according to the Riss calculus.
Operators $U$ and $S$ commute if and only if the operator $A$ is J-normal.
Moreover, the operator $A$ admits a representation
\begin{equation}
\label{f5_38}
A = U_1 S_1,
\end{equation}
where $U_1$ is a bounded J-unitary operator in $H$, and
$S_1 = \sqrt{A^T A}$ is a bounded J-self-adjoint operator in $H$.
Operators $U_1$ and $S_1$ commute if and only if $A$ is J-normal.

In particular, representations~(\ref{f5_36}) and~(\ref{f5_38}) are true for operators
\begin{equation}
\label{f5_40}
A = E + K,
\end{equation}
where $K$ is a compact operator in $H$, $\| K \| < 1$.
\end{thm}
{\bf Proof. }
Consider an operator $A$ such as in the statement of the Theorem.
We set by definition
\begin{equation}
\label{f5_40_1}
S = \sqrt{A A^T} = \int_{\Gamma} \sqrt{\lambda} R_\lambda (AA^T) d\lambda.
\end{equation}
A contour $\Gamma$ is constructed in the following way.
Let $T_R = \{ z\in{\mathbb C}:\ |z|=R \}$ be a circle, which contains $\sigma(AA^T)$ inside, $R>0$.
Let $d>0$ be a distance between a closed set
$\sigma(AA^T)$ and a segment $[0,Re^{i\varphi}]$, where $\varphi$ is from the statement of the  Theorem.
Consider parallel segments on the distance $\frac{d}{2}$ of the above segment,
join them by a half of a circle in a neighborhood of zero and completing the contour
with a part of big crcle $T_R$, it is not hard
to construct a contour $\Gamma$,
which contains the spectrum of the  operator $AA^T$ inside,
but do not contain the ray $L_\varphi$ inside.
We choose and fix an arbitrary analytic branch of the root in ${\mathbb C} \backslash L_\varphi$.

A bounded operator $B:=AA^T$ is J-self-adjoint, as it was noticed above.
Consequently, its resolvent is also a J-self-adjoint operator.
In fact, we can write
$$ R_\lambda^* (B) = ((B-\lambda E)^{-1})^* = (B^* - \overline{\lambda} E)^{-1} = (\widetilde{B} - \overline{\lambda} E)^{-1} = $$
$$= (J(B-\lambda E)J)^{-1} = J (B-\lambda E)^{-1} J = J R_\lambda (B) J,\quad \lambda\in\rho(B). $$
The operator $S$ is J-self-adjoint, as a limit of J-self-adjoint integral sums.
Moreover, there exists an inverse operator $S^{-1}$, which is also J-self-adjoint.
Set
\begin{equation}
\label{f5_41}
U = S^{-1} A,
\end{equation}
and notice that $U^{-1} = A^{-1} S$ (recall that $0\notin \sigma(A)$).
Then
$$ U \widetilde{U^*} = S^{-1} A \widetilde{A^*} \widetilde{(S^{-1})^*} = S^{-1} S^2 S^{-1} = E. $$
Multiplying the latter equality from the left side by $U^{-1}$ we get
$$\widetilde{U^*} = U^{-1}. $$
Thus, the operator $U$ is J-unitary.

Suppose now that in representation~(\ref{f5_36}) operators $U$ and $S$ commute. Then
$$ A A^T = SU \widetilde{(U^*)} \widetilde{(S^*)} = S^2, $$
$$ A^T A = \widetilde{(U^*)} S SU = S^2. $$
Conversely, if
operators $A$ and $A^T$ commute, then using last relations
(without the latter equality) we write:
$$ S^2 = \widetilde{(U^*)} S^2 U = U^{-1} S^2 U, $$
\begin{equation}
\label{f5_42}
U S^2 = S^2 U.
\end{equation}
Since $U$ commutes with $S^2$, then it commutes with an arbitrary function of this operator.
In particular, $U$ commutes with $S$.

We shall now establish a possibility of resolution~(\ref{f5_38}) for the operator $A$.
First of all we notice that
for an arbitrary linear bounded operator $D$ in $H$ we can write
$$ J R_\lambda^* (D) J = J (D^* - \overline{\lambda} E)^{-1} J = ( JD^* J - \lambda E)^{-1} =
R_\lambda (D^T),\quad \lambda\in \rho(D). $$
Therefore
\begin{equation}
\label{f5_42_1}
\rho(D) = \rho(D^T),
\end{equation}
for an arbitrary linear bounded operator $D$ in $H$.
Using this equality for operators $A$ and $AA^T$ we conclude that
$0\notin \sigma(A^T)$ and the ray $L_\varphi$ does not intersect with the spectrum of the  operator $A^T A$.
Applying the proven part of the Theorem with the operator $A^T$, we shall get a resolution
$A^T = SU$, where $S = \sqrt{A^T A}$ is a bounded J-self-adjoint operator, $U$ is
a bounded J-unitary operator.
Therefore
$$A = \widetilde{U^*} \widetilde{S^*} = U^{-1}S, $$
and it remains to notice that $U^{-1}$ is
a bounded J-unitary operator.

If the operator $A$ has the form~(\ref{f5_40}), then
$0 \notin \sigma (A)$ and
\begin{equation}
\label{f5_43}
AA^T = (E+K) J (E+K^*) J = E + C,
\end{equation}
where $C := K + JK^*J + KJK^*J$.
Notice that the operator $C$ is compact as a sum of compact operators.
The operator $J (E+K^*)^{-1}  J (E+K)^{-1}$, as it is easy to see, is the inverse operator for the operator
$AA^T$. Therefore $0 \notin \sigma(AA^T)$.
Since the spectrum of a compact operator $C$ is discrete,
having a unique point of concentration $0$, one can find a ray which is required
in the statement of the Theorem.
$\Box$
\section{Matrix representations of J-symmetric and J-skew-symmetric operators.}
We shall now turn to a study of matrix representations of J-symmetric and J-skew-symmetric
operators. Properties which are analogous to the properties of symmetric operators are valid here.
Let $J$ be a conjugation in a Hilbert space $H$ and
$\mathcal F = \{ f_k \}_{k\in{\mathbb Z}_+}$ be an orthonormal basis in $H$,
which corresponds to $J$.
Let $A$ be a linear operator in $H$, which is
J-symmetric (J-skew-symmetric) and such that
$\mathcal F \subset D(A)$.

Define a matrix of the operator $A$ in the basis $\mathcal F$:
$A_M := (a_{i,j})_{i,j\in{\mathbb Z}_+}$, $a_{i,j}=(Af_j,f_i)$.
It is not hard to verify that this matrix is complex symmetric (skew-symmetric)
in the case of J-symmetric (respectively J-skew-symmetric) operator $A$.
Notice that the columns of this matrix are square summable, i.e. belong to
$l^2$.

It is known that for an arbitrary linear operator $A$ in a Hilbert space
$H$, in the case when the set $D(A)\cap D(A^*)$ is densee in $H$,
the action of the operator $A$ is given by a matrix multirlication \cite{Cit_13000_StoneLO}.
In particular, it is true for symmetric operators.
As far as we  know, for other classes of operators a possibility to describe the action of the operator
as a matrix multiplication was not established earlier.
This property possess J-symmetric and J-skew-symmetric operators, as it shows
the following theorem.
\begin{thm}
\label{t3_4}
Let $J$ be a conjugation in a Hilbert space $H$ and
$\mathcal F = \{ f_k \}_{k\in{\mathbb Z}_+}$ be an orthonormal basis
in $H$, which corresponds to $J$.
Let $A$ be a linear operator in $H$, which is
J-symmetric (J-skew-symmetric) and such that
$\mathcal F \subset D(A)$.
Let $A_M = (a_{i,j})_{i,j\in{\mathbb Z}_+}$ be a matrix of the operator $A$ in the basis $\mathcal F$.
Then
\begin{equation}
\label{f3_26}
A g = \sum_{i=0}^\infty y_i f_i,\quad y_i=\sum_{k=0}^\infty a_{i,k} g_k,\
g=\sum_{k=0}^\infty g_k f_k \in D(A).
\end{equation}
\end{thm}
{\bf Proof. }
Let us verify the validity of the statement of the Theorem for
J-skew-symmetric operator.
For the case of J-symmetric operator the proof is analogous.
Choose an arbitrary element $g=\sum_{k=0}^\infty g_k f_k \in D(A)$.
Using that the matrix $A_M$ is skew-symmetric and using relation~(\ref{f1_4}) we write
$$y_i = (Ag, J f_i) = -(A f_i, J g) =
- (\sum_{k=0}^\infty (Af_i,f_k) f_k, \sum_{l=0}^\infty \overline{g_l} f_l) = $$
$$= - \sum_{k=0}^\infty (Af_i,f_k) g_k = - \sum_{k=0}^\infty a_{k,i} g_k =
\sum_{k=0}^\infty a_{i,k} g_k. $$
$\Box$

Let us find out, how strong the matrix $A_M$ of the operator $A$ (considered above)
determines the operator $A$.
Since J-symmetric and J-skew-symmetric operators admit closures, which are
also J-symmetric (respectively J-skew-symmetric) operators, we shall
already suppose that the operator $A$ is closed.
By the matrix $A_M$ one can define, as a matrix multiplication, an operator $T$
on $L:= \mathop{\rm Lin}\nolimits \mathcal F$.
It is easy to check that this operator is J-symmetric (J-skew-symmetric) in the case of
J-symmetric (respectively J-skew-symmetric) operator $A$.
This operator admits a closure $\overline{T}$, which is also
a J-symmetric (J-skew-symmetric) operator.
If $A=\overline{T}$, then the basis $\mathcal F$ we shall call {\bf a basis
of the matrix representation} of the operator $A$.

A question appears:
If for every complex symmetric (skew-symmetric) semi-infinite matrix $B$
with square summable columns there exists a J-symmetric (respectively J-skew-symmetric)
operator $A$ such that the matrix $B$ will be
a matrix of the operator in a corresponding to $J$ basis $\mathcal F$, and
also $\mathcal F$ will be
a basis of the matrix representation for the operator $A$?
The answer on this question is affirmative.
\begin{thm}
\label{t3_5}
Let an arbitrary complex semi-infinite symmetric (skew-symmetric) matrix
$M=(m_{i,j})_{i,j\in{\mathbb Z}_+}$ with columns in $l^2$ is given.
Then there exist a Hilbert space $H$,
a conjugation $J$ in $H$,
a J-symmetric (respectively J-skew-symmetric) operator in
$H$, a corresponding to $J$ orthonormal basis $\mathcal F$ in $H$,  $\mathcal F\subset D(A)$,
such that the matrix
$M$ is a matrix of the operator $A$ in the basis $\mathcal F$
and $\mathcal F$ is a basis of the matrix representation for $A$.
\end{thm}
{\it Proof. } For an arbitrary complex semi-infinite symmetric (skew-symmetric) matrix
$M$ with columns in $l^2$ it is enough to choose an arbitrary Hilbert space
$H$, an arbitrary orthonormal basis  $\mathcal F$ in it and to define a conjugation in
$H$ by formula~(\ref{f1_1}). Then, by using the described  above procedure, one constructs an
operator $\overline{T}$, which is the required operator.
$\Box$

Notice that, if $\mathcal F$ is a basis of the matrix representation for a closed
J-symmetric (J-skew-symmetric) operator $A$,
then $\mathcal F$ {\it will be a basis of the matrix representation for the J-adjoint operator
$\widetilde A = JAJ$, as well}.
In fact, the operator $\widetilde A$ is J-symmetric (respectively J-skew-symmetric)
by Proposition~\ref{p3_4}.
From the continuity of the operator $J$ it follows that $\widetilde A$ is closed.
Then if we choose an arbitrary element
$x\in D(\widetilde A)$, then $Jx\in D(A)$ and there exists  a sequence
$\widehat x_n\in L:=\mathop{\rm Lin}\nolimits\{ f_k \}_{k\in{\mathbb Z}_+},\ n\in{\mathbb Z}_+$:
$\widehat x_n \rightarrow Jx$, $A\widehat x_n \rightarrow AJx,\ n\rightarrow\infty$.
But then we have $J\widehat x_n\in L$, $J\widehat x_n \rightarrow x$,
$J A \widehat x_n = \widetilde A J \widehat x_n \rightarrow JAJx =
\widetilde A x,\ n\rightarrow\infty$.

The following theorem is true:

\begin{thm}
\label{t3_6}
Let $J$ be a conjugation in a Hilbert space $H$
and $\mathcal F = \{ f_k \}_{k\in{\mathbb Z}_+}$ is a corresponding to $J$ orthonormal basis in $H$.
Suppose that $A$ is a closed J-symmetric (J-skew-symmetric)
operator in $H$, $\mathcal F\subset D(A),$ and $\mathcal F$ is
a basis of the matrix representation for the operator
$A$. Let $a_{i,j} = (A f_j, f_i),\ i,j\in{\mathbb Z}_+$.
Define an operator $B$ in the following way:
\begin{equation}
\label{f3_27}
B g = \sum_{i=0}^\infty y_i f_i,\qquad y_i = \sum_{k=0}^\infty a_{i,k} g_k,\quad
g = \sum_{k=0}^\infty g_k f_k\in D_B,
\end{equation}
on a set $D_B = \{ g = \sum_{k=0}^\infty g_k f_k \in H:\ \sum_{i=0}^\infty
| \sum_{k=0}^\infty a_{i,k} g_k |^2 < \infty \}$.

\noindent
Then $A \subseteq A^T = B$ (respectively $A \subseteq - A^T = B$).

Without conditions that $A$ is closed and $\mathcal F$ is a basis of the matrix
representation for $A$, one can only state that
$A \subseteq A^T \subseteq B$ (respectively
$A \subseteq - A^T \subseteq B$).
\end{thm}
{\bf Proof. }
The proof will be given in the case of J-skew-self-adjoint operator $A$.
The case of J-symmetric operator is considered analogously.
We first show that $-A^T = -(\widetilde A)^*\subseteq B$.
Choose an arbitrary $g\in D(-(\widetilde A)^*)$ and set $-(\widetilde A)^* g = g^*$.
Let $g = \sum_{k=0}^\infty g_k f_k$, $g^* = \sum_{i=0}^\infty \widehat y_i f_i$.
We can write
$$\widehat y_i = (g^*,f_i) = (-(\widetilde A)^* g,f_i) = -(g,\widetilde A f_i) =
- (\sum_{k=0}^\infty g_k f_k, \sum_{j=0}^\infty (\widetilde A f_i,f_j) f_j) = $$
$$=- \sum_{k=0}^\infty g_k \overline{ (\widetilde A f_i,f_k) } =
- \sum_{k=0}^\infty a_{k,i} g_k = \sum_{k=0}^\infty a_{i,k} g_k,\quad i\in{\mathbb Z}_+. $$
Therefore $\sum_{i=0}^\infty
| \sum_{k=0}^\infty a_{i,k} g_k |^2 < \infty$ and, hence, we get $g\in D_B$.
Also we have $-(\widetilde A)^* g = g^* = Bg$.
Thus, we obtain an inclusion $-(\widetilde A)^*\subseteq B$.
Here we did not use that $A$ is closed and that
$\mathcal F$ is a basis of the matrix representation for $A$.
The inclusion $A \subseteq -(\widetilde A)^*$ is obvious.

Let us prove the inclusion $B\subseteq -A^T$.
As it was shown above, the operator $\widetilde A$ is closed
and $\mathcal F$ is a basis of the matrix representation for $\widetilde A$, as well.
Choose an arbitrary $g\in D_B,\
g = \sum_{k=0}^\infty g_k f_k$. Using the fact that the matrix of the operator $A$ is skew-symmetric, we write
$$(\widetilde A f_i, g) = ( \sum_{j=0}^\infty (\widetilde A f_i,f_j) f_j,
\sum_{k=0}^\infty g_k f_k) =
\sum_{k=0}^\infty (\widetilde A f_i,f_k) \overline{g_k} = $$
$$= \sum_{k=0}^\infty \overline{a_{k,i} g_k} = -
\sum_{k=0}^\infty \overline{a_{i,k} g_k} =
- \overline{\sum_{k=0}^\infty a_{i,k} g_k} = -\overline{y_i},\quad i\in{\mathbb Z}_+; $$
$$(Bg,f_i) = y_i,\quad i\in{\mathbb Z}_+. $$
Therefore
$$-(\widetilde A f_i, g) = \overline{(Bg,f_i)} = (f_i,Bg), $$
and
$$ -(\widetilde A f, g) = (f,Bg),\qquad f\in\mathop{\rm Lin}\nolimits\{ f_k \}_{k\in{\mathbb Z}_+} =: L. $$
For an arbitrary $f\in D(\widetilde A)$ there exists aa sequence
$\{ f^k \}_{k\in{\mathbb Z}_+}, f^k\in L$: $f^k\rightarrow f,\
\widetilde A f^k\rightarrow \widetilde A f$, as $k\rightarrow\infty$.
Passing to the limit as $k\rightarrow\infty$ in the equality
$$-(\widetilde A f^k, g) = (f^k,Bg)$$
and using the continuity of the scalar product, we obtain
$$ -(\widetilde A f, g) = (f,Bg),\qquad f\in D(\widetilde A). $$
Thus, we have $g\in D((\widetilde A)^*)$ è $(\widetilde A)^* g = -Bg$.
Therefore we get an inclusion $B\subseteq -(\widetilde A)^*$.
$\Box$

Let $J$ be a conjugation in  a Hilbert space $H$
and $\mathcal F = \{ f_k \}_{k\in{\mathbb Z}_+}$ be a corresponding to $J$ orthonormal basis in $H$.
Let $A$ be a  closed J-symmetric (J-skew-symmetric)
operator in $H$ and $\mathcal F\subset D(A)$. Set $a_{i,j} = (A f_j, f_i),\ i,j\in{\mathbb Z}_+,$ and
define an operator $B$ by formula~(\ref{f3_27}).
Is the operator $B$ J-symmetric (J-skew-symmetric)?
We first notice that the domain of an operator $\widetilde B = JBJ$ is a set
$$D(\widetilde B) =
\{ h=\sum_{k=0}^\infty h_k f_k \in H:\  \sum_{i=0}^\infty |\sum_{k=0}^\infty
\overline{a_{i,k}} h_k |^2 < \infty \}. $$
If $h = \sum_{k=0}^\infty h_k f_k \in D(\widetilde B)$, then
$$\widetilde B h =
\sum_{i=0}^\infty (\sum_{k=0}^\infty \overline{a_{i,k}} h_k) f_i. $$
Choose an arbitrary elements $g = \sum_{k=0}^\infty g_k f_k \in D_B$
è $h = \sum_{k=0}^\infty h_k f_k \in D(\widetilde B)$.
Using relations~(\ref{f3_7}),(\ref{f3_8}) it is easy to check that the operator $B$ is
J-symmetric (J-skew-symmetric), if the following equalities are true
(for all $g\in D_B, h\in D(\widetilde{B})$)
$$\sum_{i=0}^\infty \sum_{k=0}^\infty a_{i,k} g_k \overline{h_i} =
\sum_{k=0}^\infty
\sum_{i=0}^\infty a_{i,k} g_k \overline{h_i}. $$
In the latter case, the last theorem can be applied with the operator $B$ to obtain
that the operator $B$ is J-self-adjoint (J-skew-self-adjoint).

A question appears about existence of a basis of the matrix representation for a
closed J-symmetric (J-skew-symmetric) operator.
For an arbitrary closed operator there exists an orthonormal basis
in which the operator is a closure of its values on the linear span of the basis
(see the proof for symmetric operators in~\cite{Cit_12000_AkhGlazmanT1TLO1977}, which
is valid in the general case, as well).
A difficulty in the case of J-symmetric (J-skew-symmetric) operators is that
this new basis can be a basis which does not correspond to the  conjugation J.
So, this question remains open.

\section{A structure of the null set.}

Consider an arbitrary Hilbert space $H$. Let
$J$ be a conjugation in $H$ and $\mathcal F = \{ f_k \}_{k\in{\mathbb Z}_+}$ be a corresponding to $J$
orthonormal basis in $H$.
Let us study the set $H_{J;0}$, which we defined above ($H_{J;0} = \{ x\in H:\ [x,x]_J = 0 \}$).
Set
\begin{equation}
\label{f4_19}
H_{R} := \{ x\in H:\ (x,f_k)\in{\mathbb R},\ k\in{\mathbb Z}_+ \}.
\end{equation}
Notice that for an arbitrary element
$x\in H$ we can write a resolution:
\begin{equation}
\label{f4_20}
x = x_R + i x_I,\qquad x_R,x_I\in H_R.
\end{equation}
Namely, if $x=\sum_{k=0}^\infty x_k f_k$, we set
$x_R := \sum_{k=0}^\infty \mathop{\rm Re}\nolimits x_k f_k$, $x_I := \sum_{k=0}^\infty \mathop{\rm Im}\nolimits x_k f_k$.
It is easy to see that  representation~(\ref{f4_20}) is unique.

Define the following vectors:
\begin{equation}
\label{f4_22}
f_{k,l}^+ :=\frac{1}{\sqrt{2}} (f_k + i f_l),\quad
f_{k,l}^- :=\frac{1}{\sqrt{2}} (f_k - i f_l),\qquad k,l\in{\mathbb Z}_+.
\end{equation}
The following theorem holds true.
\begin{thm}
\label{t4_1}
Let $H$ be a Hilbert space and $J$ be a conjugation in $H$.
Let $\mathcal F = \{ f_k \}_{k=0}^\infty$ be a correesponding to $J$ orthonormal basis in $H$.
The set $H_{J;0}$ has the following properties:

\noindent
{\rm 1. } The set $H_{J;0}$ is closed;

\noindent
{\rm 2. } $x\in H_{J;0}$ $\Rightarrow$
$J x\in H_{J;0},\ \alpha x\in H_{J;0},\ \alpha\in{\mathbb C}$;

\noindent
{\rm 3. } $x,y\in H_{J;0}:\ x\perp_J y$ $\Rightarrow$
$\alpha x + \beta y \in H_{J;0},\ \alpha,\beta\in{\mathbb C}$;

\noindent
{\rm 4. }
$H_{J;0} = \{ x\in H:\ x=x_R + i x_I,\ x_R,x_I\in H_R,\
\| x_R \| = \| x_I \|,\ (x_R,x_I)=0 \}$;

\noindent
{\rm 5. } The set $H_{J;0}$ has no inner points;

\noindent
{\rm 6. } $\mathop{\rm span}\nolimits H_{J;0} = H$;

\noindent
{\rm 7. } A set $\{ f_{2k,2k+1}^+, f_{2k,2k+1}^- \}_{k\in{\mathbb Z}_+}$ is an orthonormal basis in
$H$ which elements belong to $H_{J;0}$.
\end{thm}
{\bf Proof. }
The 1-st statement of the Theorem follows from the continuity of the operator $J$ and from the
continuity of the scalar product in $H$.

\noindent
The second and third statements follows from the linearity of the  J-form and from the
properties of the conjugation $J$.

\noindent
The 4-th statement is directly verified.

\noindent
Suppose that the set $H_{J;0}$ has an inner point $x_0$ such that
\begin{equation}
\label{f4_23}
x\in H,\ \| x-x_0 \| < \varepsilon\ \Rightarrow\ x\in H_{J;0},
\end{equation}
for a number $\varepsilon>0$.
Let us write for $x_0$ the resolution~(\ref{f4_20}):
\begin{equation}
\label{f4_24}
x_0 = x_{0,R} + i x_{0,I},\qquad x_{0,R},x_{0,I}\in H_R.
\end{equation}
Suppose first that $x_{0,I}\not=0$. Set
\begin{equation}
\label{f4_25}
x_\varepsilon := x_0 + i \frac{\varepsilon}{2\| x_{0,I} \|} x_{0,I} =
x_{0,R} + i x_{0,I} \left( 1+\frac{\varepsilon}{2\| x_{0,I} \|} \right).
\end{equation}
Notice that $\| x_\varepsilon - x_0 \| = \frac{\varepsilon}{2} <\varepsilon$,
and, thus, by~(\ref{f4_23}), we obtain that $x_\varepsilon \in H_{J;0}$.
Using the proven fourth statement of the Theorem for points $x_0$ and $x_\varepsilon$, we  get
\begin{equation}
\label{f4_25_1}
\| x_{0,R} \| = \| x_{0,I} \|,
\end{equation}
and
$$\| x_{0,R} \| =
\| x_{0,I} \left( 1+\frac{\varepsilon}{2\| x_{0,I} \|} \right) \| =
\| x_{0,I} \| + \frac{\varepsilon}{2} > \| x_{0,I} \|, $$
respectively.
The obtained contradiction proves statement~5 for the case $x_{0,I}\not=0$.

If $x_{0,I}=0$, then by the fourth statement of the Theorem the relation~(\ref{f4_25_1})
is true and therefore $x_0=0$. But if zero is an inner point of the set
$H_{J;0}$, then by the proven second statement of the Theorem we get
$H_{J;0}=H$. But it is a nonsense, since, for example, elements of the basis $\mathcal F$
do not belong to the set $H_{J;0}$.

Let us prove the seventh statement of the Theorem.
Using orthonormality of elements
$f_k,\ k\in{\mathbb Z}_+,$ it is directly verified that
elements of the set
$\{ f_{2k,2k+1}^+, f_{2k,2k+1}^- \}_{k\in{\mathbb Z}_+}$, are orthonormal.
Notice that
\begin{equation}
\label{f4_27}
f_{2k} = \frac{1}{\sqrt{2}} (f_{2k,2k+1}^+ + f_{2k,2k+1}^-),\quad
f_{2k+1} = \frac{1}{\sqrt{2} i} (f_{2k,2k+1}^+ - f_{2k,2k+1}^-),\qquad k\in{\mathbb Z}_+.
\end{equation}
Therefore $\mathop{\rm span}\nolimits \{ f_{2k,2k+1}^+, f_{2k,2k+1}^- \}_{k\in{\mathbb Z}_+} = H$
and a set $\{ f_{2k,2k+1}^+, f_{2k,2k+1}^- \}_{k\in{\mathbb Z}_+}$ is an orthonormal basis in $H$.
It remains to notice that
$$[f_{2k,2k+1}^\pm,f_{2k,2k+1}^\pm ]_J =
\frac{1}{2} [f_{2k}\pm i f_{2k+1}, f_{2k}\pm i f_{2k+1}]_J =0, $$
and therefore $f_{2k,2k+1}^\pm \in H_{J;0},\ k\in{\mathbb Z}_+$.

The sixth statement of the Theorem follows from the proven seventh statement.
$\Box$

\begin{center}
\bf
On a J-polar decomposition of a bounded operator and matrix representations of J-symmetric,
J-skew-symmetric operators.
\end{center}
\begin{center}
\bf
S.M. Zagorodnyuk
\end{center}

\noindent
In this work a possibility of a decomposition of a  bounded operator which acts in a Hilbert space $H$
as a product of a J-unitary and a J-self-adjoint operators is studied, $J$ is a conjugation
(an antilinear involution).
Decompositions of J-unitary and unitary operators which are analogous to decompositions
in the finite-dimensional case are obtained. A possibility of a matrix representation for
J-symmetric, J-skew-symmetric operators is studied. Also, some simple properties
of J-symmetric, J-antisymmetric, J-isometric operators are obtained, a structure of a null set
for a J-form is studied.

\noindent
Key words and phrases: polar decomposition, matrix of an operator, conjugation, J-symmetric operator.

\noindent
MSC 2000: 47B99

}

\begin{thebibliography}{9}
\bibitem{Cit_1000_Gantmacher} F.R. Gantmaher, {\it Theory of matrices.} "Nauka", Moscow, 1967 (Russian).
\bibitem{Cit_2000_GlazmanStDAN} I.M. Glazman, {\it On an analog of the theory of extensions of
Hermitian operators and a non-symmetric one-dimensional
boundary problem on a semi-axis}// DAN SSSR {\bf 115} (1957), ¹2, 214-216 (Russian).
\bibitem{Cit_3000_GlazmanKniga} I.M. Glazman, {\it Direect methods
of qualititative spectral analysis of singular differential operators.}
Gos. izdat. fiz.-mat. liter., Moscow, 1963 (Russian).
\bibitem{Cit_4000_GarciaPutinarCS2} S.R. Garcia, M. Putinar, {\it Complex symmetric operators and applications II}//
Transactions of the AMS {\bf 359} (2007), 8, 3913-3931.
\bibitem{Cit_5000_AsadiLutsenko} Sh. Asadi, I.E. Lutsenko, {\it Skew-unitary transformations
of linear operators}//
Vestnik Kharkovskogo universiteta, Mehanika i matematika {\bf 37} (1972), ¹83, 21-27 (Russian).
\bibitem{Cit_6000_Kallinina1981} T.B. Kalinina, {\it
On extensions of an operator in a Hilbert space with a skew-unitary transformation}//
Funkcionalniy analiz (Ulyanovsk) {\bf 17} (1981), 68-75 (Russian).
\bibitem{Cit_7000_Kallinina1982} T.B. Kalinina, {\it One extension
of an operator in  a Hilbert space with a skew-unitary transformation}// Funkcionalniy analiz
(Ulyanovsk) {\bf 18} (1982), 63-71 (Russian).
\bibitem{Cit_8000_Kallinina1983} T.B. Kalinina, {\it Generalized resolvents of an operator
which is skew-symmetric with respect to an antilinear
transformation of a Hilbert space}// Funkcionalniy analiz (Ulyanovsk)
 {\bf 20} (1983), 60-72 (Russian).
\bibitem{Cit_9000_GohbergKreinTVO1967} I.T. Gohberg, M.G. Krein,
{\it Theory of Volterra operators in a Hilbert space and its applications.}
"Nauka", Moscow, 1967 (Russian).
\bibitem{Cit_10000_Kamerina1975} L.A. Kamerina, {\it Unitary equivalence of operators of the class
$K_y$}// Funkcionalniy analiz (Ulyanovsk) {\bf 5} (1975), 72-78 (Russian).
\bibitem{Cit_11000_Kamerina1987} L.A. Kamerina, {\it Quasi-unitary equivalence of operators
in spaces with an involution}// Funkcionalniy analiz (Ulyanovsk) {\bf 27} (1987), 72-78 (Russian).
\bibitem{Cit_12000_AkhGlazmanT1TLO1977} N.I. Akhiezer, I.M. Glazman, {\it Theory of linear
operators in a Hilbert space. Vol.1.} Izdat-vo pri KhGU izd. obyed. "Vysha shkola", Kharkov,
1977 (Russian).
\bibitem{Cit_13000_StoneLO} M.H. Stone, {\it Linear transformations in Hilbert space and
their applications to analysis.} AMS Colloquium Publications, Vol. 15, Providence,
Rhode Island, 1932.
\bibitem{Cit_14000_Mironov1987} B.G. Mironov, {\it On the theory of J-symmetric operators with
non-dense domain}// Funkcionalniy analiz (Ulyanovsk) {\bf 27} (1987), 128-133 (Russian).
\bibitem{Cit_15000_KrShmJpolarn} M.G. Krein, Yu.L. Shmulyan, {\it J-polar representation of plus-operators}//
Matem. issledovaniya (Kishinev), Vol.1 {\bf 2} (1966), 172-210 (Russian).
\bibitem{Cit_16000_Li1974} V.P. Li, {\it On the theory of J-symmetric operators}// Funkcionalniy analiz (Ulyanovsk)
 {\bf 3} (1974), 84-91 (Russian).
\end{thebibliography}
\end{document}